\documentclass[a4paper,11pt]{article}

\usepackage{amsthm}

\newcommand{\noi}{{\noindent}}
\newtheorem{thm}{Theorem}[section]
\newtheorem{prop}[thm]{Proposition}
\newtheorem{lem}[thm]{Lemma}

\newtheorem{coro}[thm]{Corollary}

\newtheorem{obs}[thm]{Observation}

\newtheorem{fact}[thm]{Fact}

\newcommand{\aac}{\mbox{\`a }}

\title{Permanental ideals of Hankel matrices}

\author{Elena Grieco\footnote{Dipartimento di Matematica Pura e Applicata --
Universit\aac degli Studi dell'Aquila -- Via Vetoio snc, 67010
COPPITO (AQ) Italy -- grieco@univaq.it.} , Anna
Guerrieri\footnote{Dipartimento di Matematica Pura e Applicata --
Universit\aac degli Studi dell'Aquila -- Via Vetoio snc, 67010
COPPITO (AQ) Italy -- guerrran@univaq.it.} and Irena
Swanson\footnote{Department of Mathematics -- Reed College -- 3203
SE Woodstock Blvd Portland, OR 97202 -- iswanson@reed.edu.}}
\date{}

\begin{document}

\maketitle

\begin{abstract}
We provide the Gr\"{o}bner basis and the primary decomposition of
the ideals generated by $2 \times 2$ permanents of Hankel
matrices.
\end{abstract}

\section{Introduction}

Cauchy \cite{C1812} and Binet \cite{B1812} introduced the concept
of permanent in $1812$ as a special type of alternating symmetric
function. The greater part of results on permanents in the
nineteenth century consists of identities involving permanents and
determinants. Only later the permanents were employed in various
fields of applied mathematics, in combinatorics, probability
theory, and invariant theory \cite{HK96}. A book that includes the
main results about the theory of permanents is \cite{MM78}.\\
In this paper we study permanental ideals of Hankel matrices. We
focus our attention on the structure of Gr\"{o}bner basis and
primary decomposition of the ideal generated by $2 \times 2$
permanents of Hankel matrices.\\
The {\it permanent} of an $(n \times n)$ matrix $M=(a_{ij})$ is
defined as
$$perm(M)=\sum_{\sigma \in S_n} a_{1\sigma(1)}a_{2\sigma(2)}\cdots
a_{n\sigma(n)}.$$ Thus, the permanent differs from the determinant
only in the lack of minus signs in the  expansion. If $M$ is a $u
\times v$ matrix with entries in a ring $R$, we denote by $P_r(M)$
the ideal of $R$ generated by the $(r \times r)$-subpermanents of
$M$. It is particularly interesting when $M$ is a matrix of linear
forms in a polynomial ring over a base field $K$. If the field $K$
has characteristic $2$, the permanental ideals equal the
determinantal ideals, which are well-understood \cite{C98},
\cite{W97}. Thus, from now on, we assume that the characteristic
is different from 2.\\
Laubenbacher and Swanson \cite{LS00} studied the case of generic
matrices finding the primary decomposition of $P_2(M)$. Recently
Kirkup \cite{K05} gives some indication on associated primes of
the ideal generated by $3 \times 3$ permanents of a generic matrix
(not complete list).\\ We analyze the case of Hankel matrices,
finding profoundly different results from \cite{LS00}. Let $K$ be
a field, $m,n,r$ positive integers, and $x_{i}$ variables over $K$
with $1 \leq i \leq m+n-1$. Let $R=K[x_i | 1 \leq i \leq m+n-1]$
be the polynomial ring over the previous variables, and let $M$ be
an $m \times n$ Hankel matrix in $R$, with $m \leq n$
\begin{equation}\label{mxn}
M =\left[
\begin{array}{cccc}
x_1 & x_2 & \cdots & x_n\\
x_2 & x_3 & \cdots & x_{n+1} \\
\vdots & \vdots & & \vdots \\
x_m & x_{m+1} & \cdots & x_{n+m-1}
\end{array}
\right].
\end{equation}
As $P_2(M)$ for a $2 \times 2$ Hankel matrix is a prime
ideal, in the rest of paper we assume that $m+n>4$.\\
In Section \ref{sec:grobnerbases} we analyze the structure of a
Gr\"{o}bner basis, in some cases reduced, of the ideal $P_2(M)$.
To simplify the work we start with some preliminary Lemmas about
the existence of particular monomials in the ideal $P_2(M)$,
(Lemmas \ref{lem:manymonomials1} and \ref{lem:manymonomials2}),
about the reduction of monomials of degree $2$ and $3$, (Lemma
\ref{lem:redmonomial}), and about the $S$-polynomial of two
permanents, (Lemma \ref{lem:bound}) with respect to the set $H$ of
all permanents. In the Theorem \ref{thm:gbasis2xn} we find the
reduced Gr\"{o}bner basis for $P_2(M)$ of a $2 \times n$ Hankel
matrix, and in the Theorem \ref{thm:gbasismxn} we provide a
Gr\"{o}bner basis for $P_2(M)$ of a $m \times n$ Hankel matrix
with $m \geq 3$ and $n \geq 5$. We are left with the cases $3
\times 3$, $3 \times 4$ and $4 \times 4$. Their Gr\"{o}bner bases
are different from the general case, we deal them respectively in
the Proposition \ref{prop:gbasis3x3}, ($3 \times 3$ matrix), and
in the Proposition \ref{prop:gbasis3x4&4x4}, ($3 \times 4$ and $4
\times 4$ matrices). We note that for a Hankel matrix, the
Gr\"{o}bner basis of $P_2(M)$ depends strongly on the shape of
matrix, whereas if $M$ is a $m \times n$ generic matrix there
exists a unique pattern for any $m,n$, see \cite{LS00}.\\ In
Section \ref{sec:minimalprimes} we show that the permanental ideal
$P_2(M)$ has exactly two minimal primes, see Proposition
\ref{prop:minimalprimes}, and that the structure of these primes
does not depend on the shape of matrix. This is exactly the
opposite of the situation described in \cite{LS00}. In Proposition
\ref{prop:minimalcomponents} we provide also the minimal
components of $P_2(M)$, which are different from the minimal
primes. This is in contrast with the results for generic matrices
in \cite{LS00}. In Section \ref{sec:primdec},
\ref{sec:embeddedcomp} and \ref{sec:specialcases} we provide the
full primary decomposition of $P_2(M)$. Again we need to separate
cases. In Section \ref{sec:primdec} we give a primary
decomposition, probably redundant, of the ideal $P_2(M)$. With
respect to \cite{LS00} the first difference is about the number
and the structure of minimal primes, the second difference is
about the embedded component. We show that, except for few cases,
there exists exactly one embedded component and the associated
prime is the maximal ideal. It is important to underline that, if
$M$ is a generic matrix then the embedded component exists if and
only if $m,n \geq 3$, whereas if $M$ is a Hankel matrix the
existence of the embedded component depends on the shape of
matrix, as we show in the Sections \ref{sec:embeddedcomp} and
\ref{sec:specialcases}.

\section{Gr\"{o}bner Bases}\label{sec:grobnerbases}

We first need to recall some results of Laubenbacher and Swanson
\cite{LS00} that are true also in the case of Hankel matrices. We
omit the proofs that are identical to those of \cite{LS00}.
\begin{lem}\label{lem:manymonomials1}
The ideal $P_2(M)$ contains all the products of three entries of
$M$, taken from three distinct columns and two distinct rows, or
from two distinct columns and three distinct rows.
\end{lem}
\begin{lem}\label{lem:manymonomials2}
If $m,n\geq 3$, $P_2(M)$ contains all monomials of the form
$x_{i_1 j_1}^{e_1}x_{i_2 j_2}^{e_2}x_{i_3 j_3}^{e_3}$ with
distinct $i_1, i_2, i_3$, distinct $j_1, j_2, j_3$, and where
$e_1, e_2, e_3$ are positive integers which sum to $4$.
\end{lem}
\begin{obs}
{\rm In our case $x_{ij}$ means $x_{i+j-1}$.}
\end{obs}
We start now the work to find a Gr\"{o}bner basis for $P_2(M)$
with respect to a lexicographic monomial order, with $x_1>x_2>
\cdots
>x_{m+n-1}$.
\begin{lem}\label{lem:redmonomial}
Let $M$ be a $m \times n$ Hankel matrix as in (\ref{mxn}) and let
$H$ be the set of all permanents of $M$. Define the lexicographic
ordering of monomials, with $x_1>x_2>\cdots>x_{m+n-1}$. Then for
$i, j, k =1, \ldots, m+n-1$
\begin{enumerate}
\item the reduction of $x_ix_j$ with respect to $H$ is
\begin{enumerate}
\item $(-1)^{j-1-\frac{i+j}{2}}x_{\frac{i+j}{2}}^2$, if $i+j \equiv 0 \mbox{ mod } 2$,
\item $(-1)^{j-1-\frac{i+j+1}{2}}x_{\frac{i+j-1}{2}}x_{\frac{i+j+1}{2}}$, if
$i+j\equiv 1 \mbox{ mod } 2$;
\end{enumerate}
\item the reduction of $x_ix_jx_k$ with respect to $H$ is
\begin{enumerate}
\item $\pm x_{\frac{i+j+k}{3}}^3$, if $i+j+k\equiv 0 \mbox{ mod } 3$,
\item $\pm x_{\frac{i+j+k-1}{3}}^2x_{\frac{i+j+k+2}{3}}$, if $i+j+k\equiv 1 \mbox{ mod } 3$,
\item $\pm x_{\frac{i+j+k-2}{3}}x_{\frac{i+j+k+1}{3}}^2$, if
$i+j+k\equiv 2 \mbox{ mod } 3$.
\end{enumerate}
\end{enumerate}
\end{lem}
\begin{proof} In general case we can assume that $i \leq j \leq k$.
\begin{enumerate}
\item It is clear that if $i=j$ then $i+j \equiv 0 \mbox{ mod } 2$
and $x_ix_j \equiv x_i^2$, and if $j=i+1$ then $i+j \equiv 1
\mbox{ mod } 2$ and $x_ix_j \equiv x_ix_{i+1}$. Suppose now that
$j > i+1$. By induction assumption we have
$$x_ix_j\equiv -x_{i+1}x_{j-1} \equiv -\left\{\begin{array}{lr}
(-1)^{j-2-\frac{i+j}{2}}x_{\frac{i+j}{2}}^2 & \mbox{if } i+j \equiv 0 \mbox{ mod} 2\\
(-1)^{j-2-\frac{i+j+1}{2}}x_{\frac{i+j-1}{2}}x_{\frac{i+j+1}{2}} & \mbox{if } i+j \equiv 1 \mbox{ mod} 2\\
\end{array}\right.$$
and so, we are done.
\item Clearly if $k=i$ then $i+j+k\equiv 0 \mbox{ mod }3$ and
$x_ix_jx_k=x_i^3$. Suppose $k=i+1$. Then we have two cases: if
$j=i$ then $i+j+k \equiv 1 \mbox{ mod }3$ and $x_i^2x_{i+1}$, and
if $j=k$ then $i+j+k \equiv 2 \mbox{ mod }3$ and $x_ix_{i+1}^2$.
Suppose now that $k>i+1$. Then
$$x_ix_jx_k \equiv -x_{i+1}x_jx_{k-1},$$
and for all $j=i, \ldots, k$ the difference between $max \{ i+1,
j, k-1 \} $ and $min \{ i+1, j, k-1 \} $ is strictly smaller than
$k-i$, so by induction assumption on $k-i$ we have
$$x_ix_jx_k \equiv -x_{i+1}x_jx_{k-1} \equiv
\left\{\begin{array}{lr} \pm x_{\frac{i+j+k}{3}}^3 & \mbox{if } i+j+k\equiv 0 \mbox{ mod } 3\\
\pm x_{\frac{i+j+k-1}{3}}^2x_{\frac{i+j+k+2}{3}} & \mbox{if }
i+j+k \equiv 1 \mbox{ mod } 3\\
\pm x_{\frac{i+j+k-2}{3}}x_{\frac{i+j+k+1}{3}}^2 & \mbox{if }
i+j+k\equiv 2 \mbox{ mod } 3
\end{array}
\right.$$
\end{enumerate}
and so, we are done.
\end{proof}
\begin{lem}\label{lem:bound}
Let $M$ be an $m \times n$ Hankel matrix as in (\ref{mxn}). Define
the lexicographic ordering of monomials, with
$x_1>x_2>\cdots>x_{m+n-1}$. Let $f$, $g$ be two permanents of $M$
with distinct leading terms. Then $S(f,g)$ either trivially
reduces to zero or it is a binomial whose terms $x_ax_bx_c$
satisfy $$7 \leq a+b+c \leq 3m+3n-7.$$
\end{lem}
\begin{proof}
By hypothesis on $f$ and $g$, their leading term can have only one
factor in common. Let $$f=x_ix_{i+s+t}+x_{i+s}x_{i+t},$$
$$g=x_jx_{j+s'+t'}+x_{j+s'}x_{j+t'}.$$ First suppose that $i=j$. Then
$$S(f,g)=x_{i+s}x_{i+t}x_{j+s'+t'}-x_{j+s'}x_{j+t'}x_{i+s+t}.$$
We can suppose also that $i+s+t<j+s'+t'$, so we have $$3 \leq
i+s+t \leq m+n-2,$$ $$4 \leq j+s'+t' \leq m+n-1,$$ $$1 \leq i,j
\leq m+n-4.$$ In this case the sum of indices $a+b+c$ is the same
for both terms. It is simple to see that $$ 1+3+4 \leq
i+(s+i+t)+(j+s'+t') \leq m+n-4+m+n-2+m+n-1.$$ So we have
$$8 \leq a+b+c \leq 3m+3n-7.$$ Now we suppose that $i=j+s'+t'$.
Then $$S(f,g)=x_jx_{i+s}x_{i+t}-x_{j+s'}x_{j+t'}x_{i+s+t}.$$
Clearly we have $$1 \leq j \leq i-s'-t' \leq (m+n-3)-1-1= m+n-5,$$
$$3 \leq i=j+s'+t' \leq m+n-3,$$ $$5=3+1+1 \leq (j+s'+t')+s+t=i+s+t \leq m+n-1.$$ Also
in this case the sum of indices $a+b+c$ is the same for both terms
of $S(f,g)$ and it satisfies the following relation:
$$1+3+5\leq (j)+(i)+(s+i+t) \leq
m+n-5+m+n-3+m+n-1.$$ So we have $$9 \leq a+b+c \leq 3m+3n-9.$$
Finally, we suppose that $i+s+t=j+s'+t'$. Then $$S(f,g)=
x_ix_{j+s'}x_{j+t'}+x_jx_{i+s}x_{i+t}.$$ We can suppose that
$i<j$. So we have $$1 \leq i \leq m+n-4,$$ $$2 \leq j \leq
m+n-3,$$ $$4=j+s'+t' \leq i+s+t \leq m+n-1.$$ Clearly, also in
this case the sum of indices $a+b+c$ is equal for both of terms of
$S(f,g)$ and it satisfies $$1+2+4\leq (i)+(j)+(s'+j+t') \leq
m+n-4+m+n-3-m+n-1.$$ So we have $$7 \leq a+b+c \leq 3m+3n-8.$$
\end{proof}
\begin{coro}\label{cor:corbound}
Under assumptions of Lemma \ref{lem:bound}, if a set $H$ contains
all the permanents of $M$, $x_3^3$, $x_4^3$, $\ldots$,
$x_{m+n-3}^3$, $x_2^2x_3$, $x_3^2x_4$, \ldots,
$x_{m+n-3}^2x_{m+n-2}$, $x_2x_3^2$, $x_3x_4^2$, $\ldots$,
$x_{m+n-3}x_{m+n-2}^2$, then $S(f,g)$ reduces to zero with respect
to $H$.
\end{coro}
\begin{proof}
It is sufficient to prove that each of the terms of $S(f,g)$
reduces to a multiple of some monomial of $H$. By Lemma
\ref{lem:bound}, it is sufficient to prove that whenever $i,j,k$
are positive integers with $7 \leq i+j+k \leq 3m+3n-7$, then
$x_ix_jx_k$ reduces to zero with respect to $H$. By Lemma
\ref{lem:redmonomial} $x_ix_jx_k$ reduces, with respect to
permanents, to
$$
\left.\begin{array}{lr} \pm x_{\frac{i+j+k}{3}}^3, & \mbox{if } i+j+k\equiv 0 \mbox{ mod } 3,\\
\pm x_{\frac{i+j+k-1}{3}}^2x_{\frac{i+j+k+2}{3}}, & \mbox{if }
i+j+k \equiv 1 \mbox{ mod } 3,\\
\pm x_{\frac{i+j+k-2}{3}}x_{\frac{i+j+k+1}{3}}^2, & \mbox{if }
i+j+k\equiv 2 \mbox{ mod } 3.
\end{array}
\right.$$ But the relation $7\leq i+j+k \leq 3m+3n-7$ implies
that, up to sign, these reduced monomials are in the set $x_3^3$,
$x_4^3$, \ldots, $x_{m+n-3}^3$, $x_2^2x_3$, $x_3^2x_4$, \ldots,
$x_{m+n-3}^2x_{m+n-2}$, $x_2x_3^2$, $x_3x_4^2$ \ldots,
$x_{m+n-3}x_{m+n-2}^2$, which are all elements of $H$.
\end{proof}
\begin{thm}\label{thm:gbasis2xn}
Let $M$ be a $2 \times n$ Hankel matrix, with $n \geq 3$, and let
$G$ be the following set of polynomials:
\begin{enumerate}
\item the permanents $x_ix_{i+t+1}+x_{i+1}x_{i+t}$, $i=1, \ldots,
n-1$, $t=1, \ldots, n-i$;
\item $x_i^2x_{i+1}$, $i=2, \ldots, n-1$;
\item $x_ix_{i+1}^2$, $i=2, \ldots, n-1$;
\item $x_i^3$, $i=3, \dots, n-1$;
\item $x_2^4, x_n^4$.
\end{enumerate}
Then $G$ is a minimal reduced Gr\"{o}bner basis for $P_2(M)$ with
respect to the lexicographic ordering of monomials, with
$x_1>x_2>\cdots>x_{n+1}$.
\end{thm}
\begin{proof}
First of all we observe that in case $n=3$ the set $4$ is empty.
Clearly $P_2(M)$ contains all elements of type $1$. By Lemma
\ref{lem:manymonomials1}, it is clear that $P_2(M)$ contains all
elements of type $2$ and $3$ of $G$. It is easy to prove that also
the monomials of type $4$ and $5$ are in $P_2(M)$. In fact, as the
Lemma \ref{lem:manymonomials1} assures that for all $i=3, \ldots,
n-1$, $x_{i-1}x_ix_{i+1}$, $x_2^2x_3$ and $x_{n-1}x_n^2$ are in
$P_2(M)$, we have
$$x_i^3=x_i(x_{i-1}x_{i+1}+x_i^2)-x_{i-1}x_ix_{i+1} \in P_2(M),$$
$$x_2^4=x_2^2(x_1x_3+x_2^2)-x_1x_2^2x_3 \in P_2(M),$$
$$x_n^4=x_n^2(x_{n-1}x_{n+1}+x_n^2)-x_{n-1}x_n^2x_{n+1} \in P_2(M).$$
It is easy to see also that $G$ is a reduced and minimal
generating set for $P_2(M)$. Thus it is sufficient to prove that
the $S$-polynomials of elements of $G$ reduce to zero with respect
to $G$. As the $S$-polynomial of two monomials reduces to zero, it
remains to show that $S(f,g)$
reduces to zero when $f$ is a permanent.\\
First consider $f, g$ two permanents
$$f=x_ix_{i+s+1}+x_{i+1}x_{i+s},$$
$$g=x_jx_{j+t+1}+x_{j+1}x_{j+t},$$
with $s \leq t$. Clearly $x_ix_{i+s+1}$ and $x_jx_{j+t+1}$ are
respectively the leading terms of $f$ and $g$. If $in(f)$ and
$in(g)$ have no factor in common then $S(f,g)$ reduces to zero. If
instead $in(f)$ and $in(g)$ have exactly one variable in common,
then $S(f,g)$ reduces to zero with respect to $G$, by
Corollary \ref{cor:corbound}.\\
Now, we see what happens when $f$ is a permanent and $g$ is a
monomial. We first consider $g$ a monomial of type $2$
$$f=x_ix_{i+s+1}+x_{i+1}x_{i+s},$$ $$g=x_j^2x_{j+1}.$$
If $i=j$, then $$S(f,g)=-x_ix_{i+1}^2x_{i+s}.$$ As $i=j=2, \ldots,
n-1$, it reduces to zero with respect to a monomial of type $3$.\\
If $i=j+1$, then $$S(f,g)=-x_{i-1}^2x_{i+1}x_{i+s}$$ reduces to
$x_{i-1}x_i^2x_{i+t}$ with respect to $x_{i-1}x_{i+1}+x_i^2$. As
$i=j+1=3, \ldots, n$, then it reduces to zero with respect to a
monomial of type $3$.\\
Suppose $i+s+1=j$. Clearly if $j=2$ then $i+s+t \neq j$. So
$i+s+t=j= 3, \ldots, n-1$,  and
$$S(f,g)=-x_{i+1}x_{i+s}x_{i+s+1}x_{i+s+2}.$$ But
$\alpha=x_{i+s}x_{i+s+2}+x_{i+s+1}^2$ is a permanent, so $S(f,g)$
reduces to $x_{i+1}x_{i+s+1}^3$ with respect to $\alpha$,
which reduces to zero with respect to a monomial of type $4$.\\
At the end we suppose that $i+s=j=2, \ldots, n-1$ and
$$S(f,g)=x_{i+1}x_{i+s}^3.$$
If $i+s=2$ then $S(f,g)=x_2^4$ which reduces to zero. If $i+s=3,
\ldots n-1$, $S(f,g)$ reduces to zero with respect to a monomial
of type $4$.\\
In the same way we can show that $S(f,g)$ reduces to zero when $g$
is a monomial of type $3$.\\
Now we suppose that $g$ is a monomial of type $4$, $g=x_j^3$.\\
If $i=j=2, \ldots, n-1$, then $S(f,g)=x_i^2x_{i+1}^2x_{i+s}$
reduces to zero with respect to a monomial of type $2$ or $3$.\\
If $i+s+1=j$, then $S(f,g)=x_{i+1}x_{i+s}x_{i+s+1}^2$. It is clear
that if $j=2$ then $i+s+1 \neq j$. So $i+s+1=j=3, \ldots, n-1$ and
$S(f,g)$ reduces to zero with respect to a monomial of type $3$.\\
At the end, we see what about $f$ permanent and $g=x_2^4$ or
$g=x_n^4$.\\
If $i\neq 2$, $S(f,x_2^4)=0$ and if $i=2$,
$S(f,x_2^4)=x_2^3x_3x_{2+t}$ which reduces to zero.\\
If $i+t+1 \neq n$, $S(f,x_n^4)=0$ and if $i+t+1=n$,
$S(f,x_n^4)=x_{i+1}x_{n-1}x_n^3$ which reduces to zero.
\end{proof}
\begin{prop}\label{prop:gbasis3x3}
Let $M$ be a $3 \times 3$ Hankel matrix and let $G$ be the
following set of polynomials:
\begin{enumerate}
\item the permanents $x_ix_{i+s+t}+x_{i+s}x_{i+t}$, $i=1,2,3$, $s,t=1,2$ with $i+s+t=3,4,5$;
\item $x_2^2x_{3}$, $x_2x_{3}^2$, $x_3^2x_4$, $x_3x_4^2$;
\item $x_2^4, x_3^4, x_4^4$.
\end{enumerate}
Then $G$ is a minimal reduced Gr\"{o}bner basis for $P_2(M)$ with
respect to the lexicographic ordering of monomials with $x_1>x_2>
\cdots >x_5$.
\end{prop}
\begin{proof}
Clearly $P_2(M)$ contains all elements of type $1$. By Lemma
\ref{lem:manymonomials1}, it is clear that $P_2(M)$ contains all
elements of type $2$ of $G$. Furthermore, by Lemma
\ref{lem:manymonomials1}, $x_2^2x_3$, $x_3^2x_5$,$x_3x_4^2$ are in
$P_2(M)$, then
$$x_2^4=x_2^2(x_1x_3+x_2^2)-x_1x_2^2x_3,$$
$$x_3^4=x_3^2(x_1x_5+x_3^2)-x_1x_3^2x_5,$$
$$x_4^4=x_4^2(x_3x_5+x_4^2)-x_3x_4^2x_5,$$
are in $P_2(M)$, so that also elements of type $3$ of $G$ are in
$P_2(M)$. It is easy to see also that $G$ is a reduced and a
minimal generating set for $P_2(M)$. Therefore it is sufficient to
prove that the $S$-polynomials of elements of $G$ reduce to zero
with respect to $G$. As the $S$-polynomial of two monomials
reduces to zero, it remains to show that $S(f,g)$ reduces to zero
when at least one of $f$ and $g$ is a permanent. If $in(f)$ and
$in(g)$ have no variables in common, then $S(f,g)$ always reduces
to zero. In particular,
$$S(x_1x_3+x_2^2, x_2x_4+x_3^2), S(x_1x_3+x_2^2, x_2x_5+x_3x_4),S(x_1x_3+x_2^2, x_2^4),$$
$$ S(x_1x_3+x_2^2, x_4^4), S(x_1x_4+x_2x_3, x_2x_5+x_3^2), S(x_1x_4+x_2x_3, x_3x_5+x_4^2),$$
$$S(x_1x_4+x_2x_3, x_2^2x_3), S(x_1x_4+x_2x_3, x_2x_3^2), S(x_1x_4+x_2x_3, x_2^4),$$
$$S(x_1x_4+x_2x_3, x_3^4), S(x_1x_5+x_3^2, x_2x_4+x_3^2),$$
$$S(x_1x_5+x_3^2, g) \mbox{ for all monomials } g \in G,$$
$$S(x_2x_4+x_3^2, x_3x_5+x_4^2), S(x_2x_4+x_3^2, x_3^4), S(x_2x_5+x_3x_4, x_3^2x_4),$$
$$S(x_2x_5+x_3x_4, x_3x_4^2), S(x_2x_5+x_3x_4, x_3^4), S(x_2x_5+x_3x_4, x_4^4),$$
$$S(x_3x_5+x_4^2, x_2^4), S(x_3x_5+x_4^2, x_4^4) \longrightarrow 0.$$
If $f$ and $g$ are permanents and $lcm(in(f), in(g))=x_ax_bx_c$,
with $a+b+c$ not a multiple of $3$, then as in the proof of
Corollary \ref{cor:corbound}, $S(f,g)$ reduces to zero. In
particular,
$$S(x_1x_3+x_2^2, x_1x_4+x_2x_3), S(x_1x_4+x_2x_3, x_1x_5+x_3^2),$$
$$S(x_1x_4+x_2x_3, x_2x_4+x_3^2), S(x_1x_5+x_3^2, x_2x_5+x_3x_4),$$
$$S(x_2x_4+x_3^2, x_2x_5+x_3x_4), S(x_2x_5+x_3x_4, x_3x_5+x_4^2)
\longrightarrow 0.$$ And here is the remainder of the reductions
of $S$-polynomials:
$$S(x_1x_3+x_2^2, x_1x_5+x_3^2) = x_2^2x_5-x_3^3 \stackrel{x_2x_5+x_3x_4}
\longrightarrow -x_3(x_2x_4+x_3^2) \longrightarrow 0,$$
$$S(x_1x_3+x_2^2, x_3x_5+x_4^2) = x_2^2x_5-x_1x_4^2 \stackrel{x_2x_5+x_3x_4}
\longrightarrow -x_4(x_1x_4+x_2x_3) \longrightarrow 0,$$
$$S(x_1x_3+x_2^2, x_2^2x_3) = x_2^4 \longrightarrow 0,$$
$$S(x_1x_3+x_2^2, x_2x_3^3) = x_2^3x_3 \longrightarrow 0,$$
$$S(x_1x_3+x_2^2, x_3^2x_4) = x_2^2x_3x_4 \longrightarrow 0,$$
$$S(x_1x_3+x_2^2, x_3x_4^2) = x_2^2x_4^4 \stackrel{x_2x_4+x_3^2} \longrightarrow x_3^4
\longrightarrow 0,$$
$$S(x_1x_3+x_2^2, x_3^4) = x_2^2x_3^3 \longrightarrow 0,$$
$$S(x_1x_4+x_2x_3, x_3^2x_4)= x_2x_3^3 \longrightarrow 0,$$
$$S(x_1x_4+x_2x_3, x_3x_4^2)= x_2x_3^2x_4 \longrightarrow 0,$$
$$S(x_1x_4+x_2x_3, x_4^4)= x_2x_3x_4^3 \longrightarrow 0,$$
$$S(x_1x_5+x_3^2, x_3x_5+x_4^2) = -x_1x_4^2+x_3^3 \stackrel{x_1x_4+x_2x_3}
\longrightarrow x_3(x_2x_4+x_3^2) \longrightarrow 0,$$
$$S(x_2x_4+x_3^2, x_2^2x_3) = x_2x_3^3 \longrightarrow 0,$$
$$S(x_2x_4+x_3^2, x_2x_3^2) = x_3^4 \longrightarrow 0,$$
$$S(x_2x_4+x_3^2, x_3^2x_4) = x_3^4 \longrightarrow 0,$$
$$S(x_2x_4+x_3^2, x_3x_4^2) = x_3^3x_4 \longrightarrow 0,$$
$$S(x_2x_4+x_3^2, x_2^4) = x_2^3x_3^2 \longrightarrow 0,$$
$$S(x_2x_4+x_3^2, x_4^4) = x_3^2x_4^3 \longrightarrow 0,$$
$$S(x_2x_5+x_3x_4, x_2^2x_3) = x_2x_3^2x_4 \longrightarrow 0,$$
$$S(x_2x_5+x_3x_4, x_2x_3^2) = x_3^3x_4 \longrightarrow 0,$$
$$S(x_2x_5+x_3x_4, x_2^4) = x_2^3x_3x_4 \longrightarrow 0,$$
$$S(x_3x_5+x_4^2, x_2^2x_3) = x_2^2x_4^4 \stackrel{x_2x_4+x_3^2}
\longrightarrow -x_2x_3^2x_4 \longrightarrow 0,$$
$$S(x_3x_5+x_4^2, x_2x_3^2) = x_2x_3x_4^2 \longrightarrow 0,$$
$$S(x_3x_5+x_4^2, x_3^3x_4) = x_3x_4^3 \longrightarrow 0,$$
$$S(x_3x_5+x_4^2, x_3x_4^2) = x_4^4 \longrightarrow 0,$$
$$S(x_3x_5+x_4^2, x_3^4) = x_3^3x_4^2 \longrightarrow 0.$$
\end{proof}
\begin{prop}\label{prop:gbasis3x4&4x4}
Let $M$ be a $m \times n$  Hankel matrix, where $(m,n)=(3,4)$ or
$(m,n)=(4,4)$ and let $G$ be the following set of polynomials:
\begin{enumerate}
\item the permanents $x_ix_{i+s+t}+x_{i+s}x_{i+t}$, $i=1, \ldots, m+n-3$,
$s=1, \ldots, m-1$, $t=1, \ldots n-1$ with $i+s+t=3, \ldots,
m+n-1$;
\item $x_2^2x_{3}$, $x_{m+n-3}x_{m+n-2}^2$;
\item $x_3^2, \ldots, x_{m+n-3}^2$;
\item $x_2^4, x_{m+n-2}^4$.
\end{enumerate}
Then $G$ is a minimal reduced Gr\"{o}bner basis for $P_2(M)$ with
respect to the lexicographic ordering of monomials, with $x_1>x_2>
\cdots >x_{m+n-1}$.
\end{prop}
\begin{proof} As the proof of Proposition \ref{prop:gbasis3x4&4x4} is similar
to that of Proposition \ref{prop:gbasis3x3}, we omit it.
\end{proof}
\begin{thm}\label{thm:gbasismxn}
Let $M$ be an $m \times n$ Hankel matrix with $m \geq 3$ and $n
\geq 5$, and let $G$ be the following set of polynomials:
\begin{enumerate}
\item the permanents $x_ix_{i+s+t}+x_{i+s}x_{i+t}$, $i=1, \ldots,
m+n-3$, $s=1, \ldots, m-1$, $t=1, \ldots, n-1$ with $i+s+t=3,
\ldots, m+n-1$;
\item $x_ix_{i+1}$, $i=3, \ldots, m+n-4$;
\item $x_2^2x_3$, $x_{m+n-3}x_{m+n-2}^2$;
\item $x_i^2$, $i=3, \dots, m+n-3$;
\item $x_2^4, x_{m+n-2}^4$.
\end{enumerate}
Then $G$ is a Gr\"{o}bner basis for $P_2(M)$ with respect to the
lexicographic ordering of monomials with
$x_1>x_2>\cdots>x_{m+n-1}$.
\end{thm}
\begin{proof}
First of all we show that $P_2(M)$ contains $G$. Clearly $P_2(M)$
contains all elements of type $1$, and by Lemma
\ref{lem:manymonomials1}, it is clear that $P_2(M)$ contains both
elements of type $3$ of $G$. It remains to prove that also the
monomials of type $2$, $4$ and $5$ are in $P_2(M)$. We can
consider the following submatrix of $M$
$$\left[
\begin{array}{ccccc}
x_{i-2} & x_{i-1} & x_{i} & x_{i+1} & x_{i+2} \\
x_{i-1} & x_{i} & x_{i+1} & x_{i+2} & x_{i+3} \\
x_{i} & x_{i+1} & x_{i+2} & x_{i+3} & x_{i+4} \\
\end{array}
\right].$$ Then, for all $i=3, \ldots, m+n-4$, we have
$$2x_ix_{i+1}=(x_{i-1}x_{i+2}+x_ix_{i+1})+(x_{i-2}x_{i+3}+x_ix_{i+1})-(x_{i-2}x_{i+3}+x_{i-1}x_{i+2}).$$
Moreover, for all $i=3, \ldots, m+n-3$ we have
$$2x_i^2=(x_{i-1}x_{i+1}+x_i^2)+(x_{i-2}x_{i+2}+x_i^2)-(x_{i-2}x_{i+2}+x_{i-1}x_{i+1}).$$
This proves that the monomials of type $2$ and $4$ are in
$P_2(M)$. Finally
$$x_2^4=x_2^2(x_1x_3+x_2^2)-x_1x_2^2x_3,$$
$$x_{m+n-2}^4=x_{m+n-2}^2(x_{m+n-3}x_{m+n-1}+x_{m+n-2}^2)-x_{m+n-3}x_{m+n-2}^2x_{m+n-1}.$$
Because $x_2^2x_3$ and $x_{m+n-2}^2x_{m+n-1}$ lie in $P_2(M)$ by
Lemma \ref{lem:manymonomials1}, we conclude that $G$ is in
$P_2(M)$. Now we show that the $S$-polynomials of elements of $G$
reduce to zero with respect to $G$. As the $S$-polynomial of two
monomials reduces to zero, it remains to show that $S(f,g)$
reduces to zero when at least one of $f$ and $g$ is a permanent.
First of all we consider the case in which both of them are
permanents. Let
$$f=x_ix_{i+s+t}+x_{i+s}x_{i+t},$$
$$g=x_jx_{j+s'+t'}+x_{j+s'}x_{j+t'},$$
with $i,j=1, \ldots, m+n-3$, $s=1, \ldots, m-1$, $s'=1, \ldots,
m-1$, $t=1, \ldots, n-1$, $t'=1, \ldots, n-1$, with $i+s+t$,
$j+s'+t' \leq m+n-1$. We can suppose also that $s \leq t$ and $s'
\leq t'$. Clearly $inf=x_ix_{i+s+t}$ and $in(g)=x_jx_{j+s'+t'}$.
If $in(f)$ and $in(g)$ have no factor in common then $S(f,g)$
reduces to zero. So we suppose that $in(f)$ and $in(g)$ have at
least one variable in common. First of all we suppose that
$in(f)=in(g)$ and $f\neq g$, so $i=j$ and $s+t=s'+t'$ but $s \neq
s'$ and $t \neq t'$. Suppose $s'<s$, then there exists $r= 1,
\ldots, s-1$ such that $s'= s-r$ and $t'=t+r$, and so we have
$$S(f,g)=x_{i+s-r}x_{i+t+r}-x_{i+s}x_{i+t}.$$
But we can consider the permanent
$\alpha=x_{i+s-r}x_{i+t+r}+x_{i+s}x_{i+t}$ and so $S(f,g)$ reduces
to $-2x_{i+s}x_{i+t}$ with respect to $\alpha$, and by Lemma
\ref{lem:redmonomial} it reduces to
$$(-1)^{\frac{t-s}{2}}2x_{i+\frac{s+t}{2}}^2 \ \ \ \mbox{ if } \ \ \ s+t \equiv 0
\mbox{ mod } 2,$$ or it reduces to
$$(-1)^{\frac{t-s-1}{2}}2x_{i+\frac{s+t-1}{2}}x_{{i+\frac{s+t+1}{2}}}
\ \ \ \mbox{ if } \ \ \ s+t \equiv 1 \mbox{ mod } 2.$$
By hypothesis on indices $s,s',t,t'$ we have
$$1 \leq i \leq m+n-5,$$
$$2 \leq s \leq m-1 \ \mbox{ and } \ 1 \leq t \leq n-1,$$
$$1 \leq s' \leq m-1 \ \mbox{ and } \ 2 \leq t' \leq n-1,$$
$$5 \leq i+s+t \leq m+n-1,$$
if $i=m+n-4$ or $i+s+t=4$, the only possible permanents with same
leading terms are equal. If $s+t \equiv 0 \mbox{ mod } 2,$ then
$$3= 1+ \frac{2+2}{2} \leq i+ \frac{s+t}{2} \leq \frac{m+n-5+
m+n-1}{2} =m+n-3,$$ so $x_{i+\frac{s+t}{2}}^2$ is a monomial of
type $4$. If $s+t \equiv 1 \mbox{ mod } 2$ then
$$3= 1+ \frac{2+3-1}{2} \leq i+ \frac{s+t-1}{2} \leq \frac{m+n-6+
m+n-1-1}{2} =m+n-4,$$ In fact $s+t$ must be odd but it cannot be
$3$ because in this case the only possible permanent with leading
terms $x_1x_3$ is $x_1x_3+x_2^2$. It follows
$x_{i+\frac{s+t-1}{2}}x_{{i+\frac{s+t+1}{2}}}$ is a monomial of
type $2$. Hence, in any case $S(f,g)$ reduces to zero with respect
to $G$.\\ By Lemma \ref{lem:bound} and Corollary
\ref{cor:corbound} the $S$-polynomial of two permanents reduces to
zero whenever the
leading terms of polynomials are different.\\
Now it remains to prove that $S(f,g)$ reduces to zero with respect
to $G$ when $f$ is a permanent and $g$ is a monomial. As before,
$$f=x_ix_{i+s+t}+x_{i+s}x_{i+t}$$ where $i=1, \ldots, m+n-3$,
$s=1, \ldots, m-1$, $t=1, \ldots, n-1$ with $i+s+t=3, \ldots
m+n-1$.\\ First of all we consider $g=x_jx_{j+1}$ a monomial of
type $2$. It is impossible that $in(f)=x_jx_{j+1}$, so the only
possibility is that $in(f)$ and $g$ have one factor in common.\\
If $i=j$ then $S(f,g)=x_{i+1}x_{i+s}x_{i+t}.$  As $3 \leq i=j \le
m+n-4$, if $s=t=1$ then $S(f,g)$ is a monomial of type $4$.
Otherwise $6 \leq i+s+t \leq m+n-1$ and by Lemma
\ref{lem:redmonomial} for $x_{i+s}x_{i+t}$, it reduces to a
multiple of
$$\pm x_{i+\frac{s+t}{2}}^2 \ \ \ \mbox{ if } \ \ \ s+t \equiv 0
\mbox{ mod } 2,$$ or it reduces to a multiple of
$$\pm x_{i+\frac{s+t-1}{2}}x_{{i+\frac{s+t+1}{2}}}
\ \ \ \mbox{ if } \ \ \ s+t \equiv 1 \mbox{ mod } 2.$$ We observe
that, if $s+t \equiv 0 \mbox{ mod } 2,$ then $$5=3+\frac{2+2}{2}
\leq i+ \frac{s+t}{2} \leq \frac{m+n-5+ m+n-1}{2} =m+n-3,$$ so
$x_{i+\frac{s+t}{2}}^2$ is a monomial of type $4$. However, if
$i=m+n-4$ and $i+s+t=m+n-1$, necessarily $i+s=m+n-2$ and
$i+t=m+n-3$ and so it is a monomial of type $2$. If $s+t \equiv 1
\mbox{ mod } 2$ then
$$4= 3+ \frac{2+1-1}{2} \leq i+ \frac{s+t-1}{2} \leq \frac{m+n-4+
m+n-1-1}{2} =m+n-3,$$ so
$x_{i+\frac{s+t-1}{2}}x_{{i+\frac{s+t+1}{2}}}$ is a monomial of
type $2$. Hence, in any case it reduces to zero with respect to
$G$. It is simple to see that also in the cases $x_i=x_{j+1}$,
$x_{i+s+t}=x_j$ or $x_{i+s+t}=x_{j+1}$ the $S$-polynomials reduce
to zero with respect to $G$.\\
Now, let $g=x_2^2x_3$. It is impossible that $x_ix_{i+s+t}=x_2x_3$
so the only possibilities are $x_i=x_2$, $x_i=x_3$ or
$x_{i+s+t}=x_3$.\\
If $i=2$ then $S(f,x_2^2x_3)=x_2x_3x_{2+s}x_{2+t}.$ Now we can
apply the Lemma \ref{lem:redmonomial}. As $m \geq 3$ and $n \geq
5$, then $$\frac{m+n}{2}< m+n-3,$$
$$3 \leq 2+\frac{s+t}{2} \leq \frac{m+n}{2}< m+n-3 \ \mbox{ if } s+t \equiv 0 \mbox{ mod }2,$$
$$3 \leq i+\frac{s+t-1}{2} \leq \frac{m+n}{2}< m+n-3<m+n-2 \ \mbox{ if } s+t \equiv 1 \mbox{ mod }2.$$
These arguments show that $S(f,g)$ reduces to zero with respect to
$G$.\\ If $i=3$, then $S(f,x_2^2x_3)=x_2^2x_{3+s}x_{3+t}.$ The
same argument shows that it reduces to zero with respect to $G$.\\
If $i+s+t=3$ then $i=s=t=1$ so $S(f,x_2^2x_3)=x_2^4$ that reduces
to zero with respect to $G$.\\
By symmetry, if $g=x_{m+n-3}x_{m+n-2}^2$ then $S(f,g)$
reduces to zero with respect to $G$.\\
If $g$ is a monomial of type $4$, then $S(f,g)=x_jx_{i+s}x_{i+t}$
not only in the case $i=j$ but also in the case $i+s+t=j$. In both
cases, by Lemma \ref{lem:redmonomial} and previous arguments, it
reduces to zero with respect to $G$.\\ It remains to consider $g$
a monomial of type $5$. If $g=x_2^4$ then, in order to have a
nontrivial $S$-polynomial, we can have only $x_i=x_2$ so $S(f,g)=
x_2^3x_{2+s}x_{2+t}$, which reduces to zero by Lemma
\ref{lem:redmonomial}. Finally, we see that if $g=x_{m+n-2}^4$
then, in order to have a nontrivial $S$-polynomial, the only
possibility is $x_{i+s+t}=x_{m+n-2}$ and so
$S(f,g)=x_{m+n-2}^3x_{m+n-2-t}x_{m+n-2-s}$ reduces to zero with
respect to $G$, by Lemma \ref{lem:redmonomial}.
\end{proof}

\section{Minimal primes and minimal components of $P_2(M)$}\label{sec:minimalprimes}

For the rest of this paper we set:
\begin{enumerate}
\item $r=m+n-2$, and to avoid the trivial case, we assume $r \geq 3$;
\item $P_1 =(x_1, \ldots, x_{r})$ and $P_2=(x_2, \ldots,
x_{r+1})$;
\item for all $r \geq 3$
$$\begin{array}{lll}
Q_1 & = & (x_1, \ldots, x_{r-3}, x_{r-2}^{2}, x_{r-2}x_{r-1},
x_{r-2}x_r, x_{r-2}x_{r+1}+x_{r-1}x_r,\\
 & & x_{r-1}^{2},x_{r-1}x_{r+1}+x_r^{2}),
\end{array}$$
$$Q_2=(x_5, \ldots, x_{r+1}, x_1x_3+x_2^{2}, x_1x_4+x_2x_3,
x_2x_4, x_3^{2}, x_3x_4, x_4^{2}).$$
\end{enumerate}
We now prove that $P_1$ and $P_2$ are the minimal primes of
$P_2(M)$ and $Q_1$ and $Q_2$ are the corresponding minimal
components.
\begin{prop}\label{prop:minimalprimes}
With $M$ an $m \times n$ Hankel matrix as in (\ref{mxn}), $m+n
\geq 5$, the prime ideals of $R$ minimal over $P_2(M)$ are
$$P_1=(x_1, x_2, \ldots, x_{m+n-2}) \mbox{ and } P_2=(x_2, x_3, \ldots,
x_{m+n-1}).$$
\end{prop}
\begin{proof}
Let $P$ be a minimal prime over $P_2(M)$. By Lemma
\ref{lem:manymonomials1}, $x_1x_3^{2} \in P_2(M)$. Then
$$x_1x_3^{2} \in P \Rightarrow x_1 \in P \mbox{ or } x_3 \in P.$$
We suppose that $x_1 \in P$. Since $x_i^{2}+ x_{i-1}x_{i+1} \in
P_2(M) \subseteq P$ for all $i=2, \ldots, m+n-2$ and $P$ is prime,
by induction on $i$ we have that $x_1$, $x_2$, \ldots, $x_{m+n-2}$
are all elements of $P$. Therefore $$P_2(M)\subseteq (x_1, \ldots,
x_{m+n-2}) \subseteq P,$$ and by minimality of $P$ we have that
one of the minimal primes is $$P =P_1 =(x_1, \ldots, x_{m+n-2}).$$
Now, if $x_1 \not \in P$ then $x_3 \in P$. We note that
$$x_1x_{i+j-1}+x_ix_j \in P_2(M) \subseteq P \ \ \mbox{ for all }
i=2, \ldots, m  \mbox{ and } j=2, \ldots, n.$$ Then the case
$i=j=2$ implies that $$x_1x_3+x_2^{2} \in P \mbox{ and } x_3 \in P
\mbox{ so that } x_2 \in P.$$ Suppose we have proved that  $x_2,
x_3, \ldots, x_r \in P$ for some $r \geq 3$. Then by choosing $i
\in \{2, \ldots, m\}, j \in \{2, \ldots, n\}$ such that
$i+j-1=r+1$, necessarily $i,j \leq r$, so that $$x_1x_{r+1}+
x_ix_j \in P_2(M) \subseteq P$$ implies that $x_{r+1} \in P$.
Therefore $$P_2(M) \subseteq (x_2, \ldots, x_{m+n-1}) \subseteq
P,$$ and by minimality of $P$ we have that the only other minimal
prime over $P_2(M)$ is $$P=P_2=(x_2, \ldots, x_{m+n-1}).$$
\end{proof}
\begin{prop}\label{prop:minimalcomponents}
The ideals $Q_1$ and $Q_2$ are respectively primary to $P_1$ and
$P_2.$
\end{prop}
\begin{proof}
By symmetry, it is sufficient to prove that $Q_1$ is primary to
$P_1$. By the structure of $Q_1$, it is sufficient to prove the
assertion only in the case $r=3$. Set $A=\{x_1^2, x_1x_2, x_1x_3,
x_2^2, x_1x_4+x_2x_3, x_2x_4+x_3^2\}$. We establish the following:
\begin{enumerate}
\item the degree lexicographic monomial ordering $x_1>x_2>x_3$ and $x_4$ treated as a
constant;
\item $P_2(M) \subseteq (A)$;
\item $(A) \subseteq P_2(M)_{(x_1, x_2, x_3)}\cap K[x_1, x_2, x_3,
x_4]$;
\item all $S$-polynomials of elements in $A$ reduce to zero with respect to
$A$;
\item the leading coefficients of elements of $A$ are elements of
$K$ (do not involve $x_4$).
\end{enumerate}
Then, by arguments in the proof of Proposition $3.6$ of Gianni,
Trager and Zacharias \cite{GTZ88}, we have that
$(A)=P_2(M)_{(x_1, x_2, x_3)}\cap K[x_1, x_2, x_3,$ $x_4]$ is the
$(x_1,x_2,x_3)$-primary component.\\
Clearly $P_2(M) \subseteq (A)$.\\
Moreover, it easy to prove that $(A) \subseteq P_2(M)_{(x_1, x_2,
x_3)}\cap K[x_1, x_2, x_3, x_4]$. In fact, by Lemma
\ref{lem:manymonomials1} the monomials $x_1x_2x_4$, $x_1x_3x_4$,
$x_2^2x_4$, are in $P_2(M)$ and so, as $x_4$ is units in $R_{(x_1,
x_2, x_3)}$, we have that $x_1x_2$, $x_1x_3$, $x_2^2$ are in
$P_2(M)_{(x_1, x_2, x_3)}$. Now, we notice that
$x_1^2x_4=x_1(x_1x_4+x_2x_3)-(x_1x_2)x_3$ is in $P_2(M)_{(x_1,
x_2, x_3)}$ and so, as $x_4$ is units in $R_{(x_1, x_2, x_3)}$
we have that $x_1^2$ is in $P_2(M)_{(x_1, x_2, x_3)}$.\\
Finally we prove that the $S$-polynomials of elements in $A$
reduce to zero with respect to $A$. As the $S$-polynomial of two
monomials reduces to zero, it is sufficient to prove that $S(f,g)$
reduces to zero when $f$ is a permanent. For example,
$$S(x_1x_4+x_2x_3, x_2x_4+x_3^2) = x_2^2x_3 - x_1x_3^2
\stackrel{x_2^2}{\longrightarrow} -x_1x_3^2
\stackrel{x_1x_3}{\longrightarrow}0.$$ The others are analogous.
\end{proof}

\section{A primary decomposition of $P_2(M)$}\label{sec:primdec}

In this section we find a redundant primary decomposition of the
ideal $P_2(M)$. We start by identifying what it will be the
embedded component in the cases in which this will be present.
\begin{prop}
The ideal $J=P_2(M)+(x_1^2,x_{r+1}^2)$ is primary to $(x_1,
\ldots, x_{r+1}).$
\end{prop}
\begin{proof}
To prove that the ideal $J$ is primary to $(x_1,x_2, \ldots,
x_{r+1})$, we compute $\sqrt{J}$. Clearly $x_1,x_{r+1}\in
\sqrt{J}$, $x_1x_3+x_2^2 \in J \subseteq \sqrt{J}$ so $x_2^2 \in
\sqrt{J}$ so $x_2 \in \sqrt{J}$. Say $x_1, \ldots, x_i \in
\sqrt{J}, i<r$. We see that $x_ix_{i+2}+x_{i+1}^2 \in J \subseteq
\sqrt{J}$ so $x_{i+1}^2 \in \sqrt{J}$ and so $x_{i+1} \in
\sqrt{J}$. So $J \subseteq (x_1, \ldots, x_{r+1}) \subseteq
\sqrt{J}$. But $(x_1, \ldots, x_{r+1})$ is maximal in $R=K[x_1,
\ldots, x_{r+1}]$, so $J$ is primary.
\end{proof}
\noi Now we recall a fact whose proof is folklore.
\begin{fact}\label{fact:gtztrick}
In an Noetherian ring $R$, for all ideals $I$ and for element $x
\not \in \sqrt I$, there exists an integer $n$ such that
$$(I:x^n)=(I:x^{n+1})$$
and then $$I=(I:x^n) \cap (I+(x^n)).$$
\end{fact}
\begin{thm}
Let $M$ be an $m \times n$ Hankel matrix in $K[x_1, x_2, \ldots,
x_{m+n-1}]$ as in (\ref{mxn}). Let $P_2(M)$ be the ideal generated
by the $2 \times 2$ permanents of $M$. Then a possibly redundant
primary decomposition of $P_2(M)$ is $P_2(M)=Q_1 \cap Q_2 \cap J$.
\end{thm}
\begin{proof}
We show that
\begin{enumerate}
\item $Q_1=(P_2(M):x_{r+1}^2)=(P_2(M):x_{r+1}^3)$,
\item $Q_2=(P_2(M)+(x_{r+1}^2)):x_{1}^2=(P_2(M)+(x_{r+1}^2)):x_{1}^3$,
\end{enumerate}
and so by Fact \ref{fact:gtztrick} we can assert that $$P_2(M)=Q_1
\cap Q_2 \cap J$$ is a primary decomposition of $P_2(M)$.
\begin{enumerate}
\item
   \begin{enumerate}
       \item First of all we see that $x_{r+1}^2Q_1 \subseteq
         P_2(M)$. As $x_{r-2}x_{r+1}+x_{r-1}x_r$, and $x_{r-1}x_{r+1}+x_r^{2}$
         are two permanents, it is clear that
         $(x_{r-2}x_{r+1}+x_{r-1}x_r)x_{r+1}^2$,
         $(x_{r-1}x_{r+1}+x_r^{2})x_{r+1}^2$ are in $P_2(M)$. By Lemma \ref{lem:manymonomials1},
         $x_{r-2}^{2}x_{r+1}^2$, $x_{r-2}x_{r-1}x_{r+1}^2$,
         $x_{r-2}x_rx_{r+1}^2$ and $x_{r-1}^{2}x_{r+1}^2$ are in
         $P_2(M)$. If $r \geq 4$ and $i=1,\ldots, r-3$ then
         $x_ix_{r+1}^2 \in P_2(M)$. In fact as $(x_ix_{r+1}+x_sx_t)$
         with $s,t$ such that $s+t=i+r+1$ is a permanent and as, by Lemma \ref{lem:manymonomials1},
         $x_sx_tx_{r+1}$ is in $P_2(M)$ we have that
         $x_ix_{r+1}^2=x_{r+1}(x_ix_{r+1}+x_sx_t)-x_sx_tx_{r+1}$
         is in $P_2(M)$.
       So we have $Q_1 \subseteq (P_2(M):x_{r+1}^2)$.
       But we know that $P_2(M) \subseteq Q_1$ and $Q_1$ is primary, so
       $Q_1 \subseteq (P_2(M):x_{r+1}^2) \subseteq (Q_1:x_{r+1}^2)= Q_1$.
       Then $(P_2(M):x_{r+1}^2)=Q_1$.
       \item For all $y\not\in \sqrt{Q_1}$ we have $Q_1 x_{r+1}^2y \subseteq P_2(M)\subseteq Q_1$
       so $Q_1=(P_2(M):x_{r+1}^2y)$ and in particular we have
       $Q_1=(P_2(M):x_{r+1}^2)=(P_2(M):x_{r+1}^3)$.
   \end{enumerate}
  \item \begin{enumerate}
\item By symmetry, $x_1^2 Q_2 \subseteq P_2(M)$. So we have $Q_2 \subseteq (P_2(M)+(x_{r+1}^2)):x_1^2$.
But $P_2(M) +(x_{r+1}^2) \subseteq Q_2$ and $Q_2$ is primary, so
$Q_2 \subseteq ((P_2(M) +(x_{r+1}^2)):x_1^2) \subseteq (Q_2:x_1^2)
=Q_2$. Then $((P_2(M) +(x_{r+1}^2)):x_1^2) =Q_2$.
\item For all $y\not\in\sqrt{Q_2}$ we have $Q_2 x_{1}^2y \subseteq P_2(M)+(x_{r+1}^2)\subseteq Q_2$
so $Q_2=((P_2(M)+(x_{r+1}^2)):x_{1}^2y)$ and in particular we have
$Q_2=((P_2(M)+(x_{r+1}^2)):x_{1}^2)=((P_2(M)+(x_{r+1}^2)):x_{1}^3)$.
  \end{enumerate}
\end{enumerate}
\end{proof}
\begin{obs}
{\rm Whereas $Q_1$ and $Q_2$ are never redundant, $J$ may be
redundant, but only in finitely many cases. We describe precisely
what happens in the next Sections.}
\end{obs}

\section{When there is an embedded component}\label{sec:embeddedcomp}

\begin{prop}
Let $M$ be an $m \times n$ Hankel matrix as in (\ref{mxn}).
\begin{enumerate}
\item If $m=2$ and $n \geq 4$, then $(x_1, x_2, \ldots, x_{n+1})$ is an associated prime
of $P_2(M)$.
\item If $m \geq 3$ and $m+n-1 \geq 9$, then $(x_1, \ldots, x_{m+n-1})$ is an associated prime of
$P_2(M)$.
\end{enumerate}
\end{prop}
\begin{proof}
By definition, an ideal $J$ is an associated prime to $P_2(M)$ if
there exists $\alpha \in K[x_1, \ldots, x_{m+n-1}]$ such that
$$J = (P_2(M): \alpha).$$
\begin{enumerate}
\item Let
$$M=\left[
\begin{array}{ccccccc}
x_1 & x_2 & x_3 & x_4 & \ldots & x_{n-1} & x_n \\
x_2 & x_3 & x_4 & x_5 & \ldots & x_n & x_{n+1}
\end{array}
\right],$$ and for all $i,j$ with $i<j, \ \ i=2, \ldots, n-2, \ \
j=4, \ldots, n,$
$$\alpha = x_ix_j.$$ Then, by Lemma \ref{lem:manymonomials1}, we have
$$(x_1, x_2, \ldots, x_{n+1}) \subseteq (P_2(M): \alpha),$$
and Theorem \ref{thm:gbasis2xn} implies that $\alpha \not\in
P_2(M)$.
\item Now let $M$ a $m \times n$ Hankel matrix as in (\ref{mxn}).
We assume $m \geq 3$ and $m+n-1 \geq 9$. For all $j=5, \ldots,
(m+n-1)-4$ set
$$\alpha =x_j.$$
Then by degree count, $\alpha \not \in P_2(M)$. So it is
sufficient to prove that $$\alpha (x_1, \ldots, x_{m+n-1})
\subseteq P_2(M).$$ By Lemma \ref{lem:redmonomial}, for all $i=1,
\ldots, m+n-1$ and $j=5, \ldots, m+n-5$ the monomials $x_ix_j$
reduce, with respect to elements of $P_2(M)$, to
$$\pm x_{\frac{i+j}{2}}^2 \ \ \mbox{ if } i+j \equiv 0 \mbox{ mod }
2,$$
$$\pm x_{\frac{i+j-1}{2}}x_{\frac{i+j+1}{2}} \ \ \mbox{ if } i+j\equiv 1 \mbox{ mod } 2.$$
Clearly, if $i+j \equiv 0 \mbox{ mod } 2$, then
$$3 \leq \frac{i+j}{2} \leq m+n-3,$$
and if $i+j \equiv 1 \mbox{ mod } 2$, then
$$3 \leq \frac{i+j-1}{2} \leq m+n-4.$$
Thus, by Theorem \ref{thm:gbasismxn}, all monomials $x_ix_j$ are
in $P_2(M)$.
\end{enumerate}
\end{proof}
\begin{prop}
Let $M$ be a $3 \times 3$ Hankel matrix. Then the ideal $(x_1,
\ldots, x_5)$ is an associated prime of $P_2(M)$.
\end{prop}
\begin{proof}
In this case
$$M= \left[ \begin{array}{ccc}
x_1 & x_2 & x_3\\
x_2 & x_3 & x_4\\
x_3 & x_4 & x_5\\
\end{array}
\right].$$ Set $\alpha=x_1x_3x_5$. Then by Lemma
\ref{lem:manymonomials1} and \ref{lem:manymonomials2},
$$(x_1, \ldots, x_5) \subseteq (P_2(M): \alpha ),$$
and by Proposition \ref{prop:gbasis3x3} $x_1x_3x_5\not\in P_2(M)$.
\end{proof}
\begin{prop}
Let $M$ be a $3 \times 4$ Hankel matrix. Then the ideal $(x_1,
\ldots, x_6)$ is an associated prime of $P_2(M)$.
\end{prop}
\begin{proof}
As in previous propositions we see that if
$$M=
\left[
\begin{array}{cccc}
x_1 & x_2 & x_3 & x_4\\
x_2 & x_3 & x_4 & x_5\\
x_3 & x_4 & x_5 & x_6\\
\end{array}
\right],$$ and if $\alpha = x_2x_5$ or $\alpha = x_3x_4$ then,
Lemma \ref{lem:manymonomials1} implies that
$$(x_1, \ldots, x_6) \subseteq (P_2(M): \alpha).$$ Proposition
\ref{prop:gbasis3x4&4x4} shows that $\alpha \not\in P_2(M)$.
\end{proof}
\begin{prop}
Let $M$ be a $4 \times 4$ Hankel matrix. Then $(x_1, \ldots, x_7)$
is an associated prime of $P_2(M)$.
\end{prop}
\begin{proof}
In this case we see that if
$$M=
\left[ \begin{array}{cccc}
x_1 & x_2 & x_3 & x_4\\
x_2 & x_3 & x_4 & x_5\\
x_3 & x_4 & x_5 & x_6\\
x_4 & x_5 & x_6 & x_7\\
\end{array} \right],$$
and if $\alpha$ is any of $x_2x_5$, $x_3x_4$, $x_3x_6$ or $x_4x_5$
then Lemma \ref{lem:manymonomials1} implies that
$$(x_1, \ldots, x_7) \subseteq (P_2(M): \alpha).$$ The Proposition
\ref{prop:gbasis3x4&4x4} shows that $\alpha \not\in P_2(M)$.
\end{proof}

\section{When there are only the minimal components}\label{sec:specialcases}

We have mentioned that in a few cases the primary decomposition of
$P_2(M)$ admits only the two minimal components. In this section
we describe all such cases. Throughout we will use the following
easy facts (for the proofs see \cite{AtM} and \cite{AL94}).
\begin{fact}\label{fact:intersection}
For all ideal $I,J,K$ then $$(I+J)\cap(I+K)= I+J\cap (I+K).$$
\end{fact}
\begin{fact}\label{fact:intergrbases}
Let $t$ be a variable over $R=K[x_1, \ldots, x_{m+n-1}]$. We
impose on $R[t]$ a monomial order such that for any $f \in R[t]
\setminus R, in(f) \not \in R$. For all pairs of ideals $I, J$ in
$R$ we can compute $I \cap J$ via Gr\"{o}bner basis. Namely, $$I
\cap J = (ItR[t]+J(t-1)R[t])\cap R.$$
\end{fact}
\begin{prop}\label{prop:primdec2x3}
Let $M$ be a $2 \times 3$ Hankel matrix, then the primary
decomposition of $M$ is $Q_1 \cap Q_2$.
\end{prop}
\begin{proof}
It is sufficient to show that $P_2(M)=Q_1 \cap Q_2$. Actually, we
can compute the intersection of the minimal components.
$$Q_1=(x_1^2, x_1x_2, x_1x_3, x_2^2, x_1x_4+x_2x_3, x_2x_4+x_3^2),$$
$$Q_2=(x_2x_4, x_3^2, x_3x_4, x_4^2, x_1x_3+x_2^2,
x_1x_4+x_2x_3).$$ By using the Fact \ref{fact:intersection}
$$\begin{array}{lcl}
Q_1 \cap Q_2 & = & (x_1x_4+x_2x_3, x_2x_4+x_3^2)+\\
 & & [(x_1^2, x_1x_2, x_1x_3, x_2^2)\cap (x_2x_4, x_3^2, x_3x_4, x_4^2, x_1x_3+x_2^2,
x_1x_4+x_2x_3)]\\
 & = &(x_1x_4+x_2x_3, x_2x_4+x_3^2,x_1x_3+x_2^2)+\\
 & & [(x_1^2, x_1x_2, x_1x_3, x_2^2) \cap (x_2x_4, x_3^2,
x_3x_4, x_4^2, x_1x_4+x_2x_3)]\\
 &=& P_2(M) + [(x_1^2, x_1x_2, x_1x_3, x_2^2) \cap (x_2x_4, x_3^2,
x_3x_4, x_4^2, x_1x_4+x_2x_3)]. \end{array}$$ By using Fact
\ref{fact:intergrbases} it is straightforward to see that the last
intersection is equal to $$(x_1^2 x_4+x_1x_2x_3, x_1x_2x_4,
x_1x_3^2, x_1x_3x_4, x_2^2x_3, x_2^2x_4).$$ Clearly
$x_1^2x_4+x_1x_2x_3$ is in $P_2(M)$. By Lemma
\ref{lem:manymonomials1} the monomials $x_1x_2x_4$, $x_1x_3^2$,
$x_1x_3x_4$, $x_2^2x_3$, $x_2^2x_4$ are in $P_2(M)$, so $$Q_1 \cap
Q_2 = P_2(M).$$
\end{proof}
\begin{prop}\label{prop:primdec3x5}
Let $M$ be a $3 \times 5$ Hankel matrix. Then the primary
decomposition of $M$ is $Q_1 \cap Q_2$.
\end{prop}
\begin{proof}
We know that
$$Q_1=(x_1, x_2, x_3, x_4^2, x_4x_5, x_4x_6, x_5^2, x_4x_7+x_5x_6,
x_5x_7+x_6^2),$$
$$Q_2=(x_5, x_6, x_7, x_2x_4, x_3^2, x_3x_4, x_4^2, x_1x_3+x_2^2,
x_1x_4+x_2x_3).$$ As $(x_4^2,x_4x_5, x_4x_6, x_5^2, x_4x_7+x_5x_6,
x_5x_7+x_6^2) \subseteq (x_5, x_6, x_7, x_4^2)$, by using the Fact
\ref{fact:intersection}, we have
$$\begin{array}{lll}
Q_1 \cap Q_2 &=& (x_4^2, x_4x_5, x_4x_6, x_5^2,
x_4x_7+x_5x_6, x_5x_7+x_6^2)+[(x_1, x_2, x_3)\cap \\
 & &(x_5, x_6, x_7, x_2x_4, x_3^2,
x_3x_4, x_4^2, x_1x_3+x_2^2, x_1x_4+x_2x_3)].
\end{array}$$
Moreover, as $(x_2x_4, x_3^2, x_3x_4, x_1x_3+x_2^2, x_1x_4+x_2x_3)
\subseteq (x_1, x_2, x_3),$ by using the Fact
\ref{fact:intersection}, we have
$$\begin{array}{lll} Q_1 \cap Q_2 &=& (x_4^2, x_4x_5,
x_4x_6, x_5^2, x_4x_7+x_5x_6, x_5x_7+x_6^2) +\\
& &  (x_2x_4, x_3^2, x_3x_4, x_1x_3+x_2^2, x_1x_4+x_2x_3)+\\
 & & [(x_1, x_2, x_3) \cap (x_4^2, x_5, x_6, x_7)].
\end{array}$$
But we know that $(x_1, x_2, x_3) \cap (x_4^2, x_5, x_6, x_7)=
(x_1, x_2, x_3)\cdot(x_4^2, x_5, x_6, x_7)$, so we have
$$\begin{array}{lll}
Q_1 \cap Q_2 &=& (x_4^2) + (x_4x_5,
x_4x_6, x_5^2, x_4x_7+x_5x_6, x_5x_7+x_6^2) + \\
& &  (x_2x_4, x_3^2, x_3x_4, x_1x_3+x_2^2, x_1x_4+x_2x_3)+(x_1x_4^2, x_2x_4^2, x_3x_4^2)+\\
 & & (x_1x_5, x_1x_6, x_1x_7, x_2x_5, x_2x_6, x_2x_7, x_3x_5,
x_3x_6, x_3x_7).
\end{array}$$
As $(x_1x_4^2, x_2x_4^2, x_3x_4^2) \subseteq (x_4^2)$, then
$$ \begin{array}{lll}
Q_1 \cap Q_2 & = & (x_1x_3+x_2^2, x_1x_4+x_2x_3, x_1x_5, x_1x_6,
x_1x_7, x_2x_4,\\
 & & x_2x_5, x_2x_6, x_2x_7, x_3^2, x_3x_4, x_3x_5, x_3x_6, x_3x_7, \\
 & & x_4^2, x_4x_5, x_4x_6, x_4x_7+x_5x_6, x_5^2, x_5x_7+x_6^2).
 \end{array}$$
Clearly $Q_1 \cap Q_2$ contains $P_2(M)$, but it is also simple to
show the opposite inclusion. It is sufficient to prove that all
monomials in $Q_1 \cap Q_2$ are in $P_2(M)$. As $x_1x_2+x_3^2$,
$x_1x_5+x_2x_4$ and $x_2x_4+x_3^2$ are in $P_2(M)$, then $x_1x_5$,
$x_2x_4$ and $x_3^2$ are in $P_2(M)$. As $x_1x_6+x_2x_5$,
$x_1x_6+x_3x_4$ and $x_2x_5+x_3x_4$ are in $P_2(M)$, then
$x_1x_6$,$x_2x_5$ and $x_3x_4$ are in $P_2(M)$. Now we see that
$x_2x_6+x_3x_5$, $x_2x_6+ x_4^2$ and $x_3x_5+x_4^2$ are elements
of $P_2(M)$, which implies that $x_2x_6$, $x_3x_5$ and $x_4^2$ are
in $P_2(M)$. So, as $x_1x_7+x_3x_5$ is in $P_2(M)$, also $x_1x_7$
is in $P_2(M)$. Another time we see that $x_2x_7+x_3x_6$,
$x_2x_7+x_4x_5$ and $x_3x_6+x_4x_5$ are in $P_2(M)$, so $x_2x_7$,
$x_3x_6$ and $x_4x_5$ are in $P_2(M)$. At the end, as
$x_3x_7+x_4x_6$, $x_3x_7+x_5^2$ and $x_4x_6+x_5^2$ are in $P_2(M)$
then $x_3x_7$, $x_4x_6$ and $x_5^2$ are in $P_2(M)$ too.
\end{proof}
\begin{prop}\label{prop:primdec3x6&4x5}
Let $M$ be a $m\times n$ Hankel matrix with $(m, n)=(3, 6), (4,
5)$. Then the primary decomposition of $M$ is $Q_1 \cap Q_2$.
\end{prop}
\begin{proof}
As in previous Proposition we show that $P_2(M)= Q_1 \cap Q_2$.
Clearly $Q_1 \cap Q_2$ contains $P_2(M)$. To show the opposite
inclusion we compute the intersection between $Q_1$ and $Q_2$. In
this case we have
$$Q_1=(x_1, x_2, x_3, x_4, x_5^2, x_5x_6,
x_5x_7, x_6^2, x_5x_8+x_6x_7, x_6x_8+x_7^2),$$
$$Q_2=(x_5, x_6, x_7, x_8, x_2x_4, x_3^2, x_3x_4, x_4^2,
x_1x_3+x_2^2, x_1x_4+x_2x_3).$$ By Fact \ref{fact:intersection} it
is simple to see that
$$\begin{array}{lll}
Q_1 \cap Q_2 & = & (x_1x_3+x_2^2, x_1x_4+x_2x_3, x_2x_4, x_3^2,
x_3x_4, x_4^2, x_5^2, x_5x_6, x_5x_7,\\
 & & x_5x_8+x_6x_7, x_6^2, x_6x_8+x_7^2) + (x_1, x_2, x_3, x_4) \cdot (x_5, x_6, x_7,
 x_8).
\end{array}$$
After computing $(x_1, x_2, x_3, x_4) \cdot (x_5, x_6, x_7, x_8)$
we have
$$\begin{array}{lll}
Q_1 \cap Q_2 & = & (x_1x_3+x_2^2, x_1x_4+x_2x_3, x_1x_5, x_1x_6,
x_1x_7, x_1x_8, x_2x_4, x_2x_5,\\
 & & x_2x_6, x_2x_7, x_2x_8, x_3^2, x_3x_4, x_3x_5, x_3x_6, x_3x_7, x_3x_8, x_4^2, x_4x_5,\\
 & & x_4x_6, x_4x_7, x_4x_8, x_5^2, x_5x_6, x_5x_7, x_5x_8+x_6x_7, x_6^2,
x_6x_8+x_7^2).
\end{array}$$
It remains to prove that all monomials in $Q_1 \cap Q_2$ are in
$P_2(M)$ too. By Lemma \ref{lem:redmonomial} it is sufficient to
prove that $x_3^2$, $x_3x_4$, $x_4^2$, $x_4x_5$, $x_5^2$,
$x_5x_6$, $x_6^2$ are in $P_2(M)$. But these are in $P_2(M)$ by
Theorem \ref{thm:gbasismxn}.
\end{proof}

\section*{Acknowledgments}
The first author thanks the Ph.D program of the University of
L'Aquila for support and the New Mexico State University for
hospitality.

\bibliographystyle{plain}
\bibliography{biblio}

\end{document}